 \documentclass[11pt]{article}
 \usepackage{mathrsfs}
 \usepackage{amsfonts}
 \usepackage[plainpages=false]{hyperref}
 \usepackage{harvard}
 \usepackage{cite}
 \usepackage{xcolor}
 \usepackage{enumerate}
  \usepackage{amsmath}

 \newtheorem{theorem}{Theorem}[section]
 \newtheorem{lemma}[theorem]{Lemma}
 \newtheorem{corol}[theorem]{Corollary}
 \newtheorem{prop}[theorem]{Proposition}
 \newtheorem{condition}[theorem]{Condition}
  \newtheorem{remark}[theorem]{Remark}
  \newtheorem{definition}[theorem]{Definition}

 \def\blemma{\begin{lemma}}\def\elemma{\end{lemma}}
 \def\bproposition{\begin{prop}}\def\eproposition{\end{prop}}
 \def\btheorem{\begin{theorem}}\def\etheorem{\end{theorem}}
 \def\bcorollary{\begin{corol}}\def\ecorollary{\end{corol}}
 \def\bcondition{\begin{condition}}\def\econdition{\end{condition}}
 \def\bremark{\begin{remark}}\def\eremark{\end{remark}}
  \def\bdefinition{\begin{definition}}\def\edefinition{\end{definition}}

 \def\blem{\begin{lemma}}\def\elem{\end{lemma}}
 \def\bprop{\begin{prop}}\def\eprop{\end{prop}}
 \def\bthm{\begin{theorem}}\def\ethm{\end{theorem}}
 \def\bcor{\begin{corol}}\def\ecor{\end{corol}}
 \def\bcon{\begin{condition}}\def\econ{\end{condition}}
 \def\brem{\begin{remark}}\def\erem{\end{remark}}

 \def\beqlb{\begin{eqnarray}}\def\eeqlb{\end{eqnarray}}
 \def\beqnn{\begin{eqnarray*}}\def\eeqnn{\end{eqnarray*}}

 \def\ar{\!\!\!&}

 \def\e{\mbox{\rm e}}

 \def\proof{\noindent{\it Proof.~~}}\def\qed{\hfill$\Box$\medskip}

 \def\<{\langle}\def\>{\rangle}

\def\az{\alpha}

\def\gz{\gamma}

\def\dz{\delta}
\def\lz{\lambda}

\def\wz{\wedge}

\def\iz{\infty}

\def\F{{\mathscr F}}
\def\R{{\mathbb R}}
\def\E{{\mathbb E}}
\def\P{{\mathbb P}}

\def\tN{{\tilde{N}}}

\def\by{\bar{y}}

\begin{document}

\centerline{\LARGE\textbf{ Explosion of continuous-state branching processes}}

\bigskip

\centerline{\LARGE\textbf{with competition in L\'{e}vy  environment}}

\vskip 0.5cm

\centerline{Rugang Ma\footnote{School of Statistics and Mathematics, Central University of Finance and Economics, Beijing, China.
Department of Mathematics and Statistics, Concordia University, Montreal, Canada.
Supported by the disciplinary funding of Central University of Finance and Economics and by NSERC (RGPIN-2021-04100).
    Email: marg@cufe.edu.cn.}
 and
 Xiaowen Zhou\footnote{Department of Mathematics and Statistics, Concordia University, Montreal, Canada. Supported by NSERC
(RGPIN-2021-04100). Email: xiaowen.zhou@concordia.ca. Corresponding author.}
 }

\vskip 0.5cm

{\narrower{\narrower

\noindent \textbf{Abstract} ~
  Using the Lyapunov criteria arguments, we find sufficient conditions on explosion/nonexplosion for continuous-state branching processes with competition in L\'evy random environment. In particular, we identify the necessary and sufficient conditions on explosion/nonexplosion when the competition function is a power function and the L\'evy measure of the associated branching mechanism is stable.

\smallskip


\smallskip

\noindent{\textbf{MR(2020) Subject Classification} }:  ~60J80;  60H20.

\par}\par}


\section{Introduction}
\setcounter{equation}{0}

Continuous-state branching processes in random environment arise as scaling limits of Bienaym\'e-Galton-Watson processes that were introduced in Smith \cite{S68} and Smith and Wilkinson \cite{SW69}; see Kurtz \cite{K78} for early work on diffusion approximations of branching processes in random environment. A recent study of the Feller branching processes in Brownian environment can be found in B$\ddot {\text o}$inghoff and Hutzenthaler \cite{BH12} where asymptotics of the survival probability is studied for different regimes. The introduction of branching processes under the continuous-state setting allows to apply the stochastic differential equations (SDEs for short) and L\'evy processes techniques in its study. We refer to Kyprianou \cite{K14} and Li \cite{L11} for comprehensive introduction on continuous-state branching processes and the associated stochastic equations.

To understand the effect of random environment on demography of the branching process,
continuous-state branching process with catastrophes was first proposed in Bansaye et al. \cite{BMS13}  as a  continuous-state branching process in L\'evy environment (CBLE for short) where the random  environment is modelled by a L\'evy process with sample paths of bounded variation.
More general CBLEs were introduced and studied in He et al. \cite{HLX18} and Palau and Pardo \cite{PP18} as  unique non-negative strong solutions to certain SDEs driven by Brownian motions and L\'evy processes associated to both the branching mechanism and the random environment. We refer to Bansaye et al. \cite{BCM19} for discussions on the convergence of discrete-state population models to  CBLEs.

The quenched Laplace transform for  the branching process in random environment can be expressed using random cumulant  semigroups conditional on the environment.
A necessary and sufficient condition  in terms of Grey's condition  was shown in \cite{HLX18} for the CBLE to become extinct.  The speed of extinction was also obtained in Bansaye et al. \cite{BPS21} for CBLE for which the L\'evy environment process fluctuates.

In another development on continuous-state branching processes, a logistic branching process was introduced in Lambert \cite{Lambert05} to incorporate  competition among individuals in the continuous-state branching process.
Foucart \cite{F19} studied the boundary behavior of continuous-state branching processes with logistic competition and obtain an integral test on explosion/nonexplosion.
A general competition mechanism was introduced in
Ba and Pardoux \cite{BP15} and Ma \cite{M15}. Under a moment condition $\int_0^\infty (z\wedge z^2)\mu (dz)<\infty $  on the L\'evy measure $\mu$ for the branching mechanism, \cite{M15} established the Lamperti transformation  between continuous-state branching processes
with competition and strong solutions of  stochastic
equations driven by L\'{e}vy processes without negative jumps; see also  Berestycki et al. \cite{BFF18}  for flows of continuous-state branching processes with competition.  We refer to Li et al. \cite{LLWZ22} for recent work on ergodic results of continuous state branching processes with immigration and competition.
The continuous state branching process with immigration and competition in a L\'{e}vy random environment was introduced in \cite{PP18}
 with its long term behaviours studied.
 The extinction and coming down from infinity behaviours  were also  studied in Leman and Pardo \cite{LP18} for CBLEs with competition.

The explosion/nonexplosion  conditions for continuous-state branching processes are well known; see Grey \cite{G74} for an integral test on Laplace exponent of the associated branching mechanism.
  An integral test on explosion/nonexplosion was further proved in Leman and Pardo \cite{LP21} for continuous-state branching process in Brownian environment with a special branching mechanism that is associated to the Laplace transform of a subordinator  and with logistic competition.
  It was also pointed out that  a continuous-state branching process in  L\'evy environment is  conservative, i.e. the explosion can not happen, if the L\'evy measure $\mu$ for the branching mechanism satisfies the moment condition; see  Lemma A.1 of \cite{BPS21}.
 On the other hand, it is known that large enough competition can prevent explosion from happening;
 see Foucart \cite{F19} and Li et al. \cite{LYZ19}.
  Some sufficient conditions of explosion were found in \cite{LYZ19}  for general continuous-state nonlinear branching processes whose competition mechanism is a general function and the L\'evy measure $\mu$ for the branching mechanism satisfies the moment condition.
To our best knowledge, the explosion/nonexplosion conditions for CBLEs  with general competition and with general L\'evy measure $\mu$  have not been studied systematically.

Integral tests on explosion/nonexplosion are not available anymore for the above-mentioned branching processes with general competitions, and  as an effective alternative, the approach of Foster-Lyapunov criteria  comes into play.
The Foster-Lyapunov criteria find successful applications in characterizing the boundary behaviours of SDEs related to the continuous-state branching processes; see Li et al. \cite{LYZ19} and Ma et al. \cite{MYZ21}. We are not aware of previous applications of Foster-Lyapunov criteria in the study of
CBLEs.

In this paper, applying the Foster-Lyapunov criteria arguments to suitable test functions  we
find sufficient conditions on explosion/nonexplosion for CBLEs with  general competition. In particular, we identify necessary and sufficient conditions on explosion/nonexplosion when the competition function is a power function and jump part of the branching mechanism is  an  $\alpha$-stable process for $\alpha\in(0,2)$, which  helps to  determine the interplay between competition and large jumps of the branching on the explosion.  As a corollary we also show that the Neveu's CBLE with competition can not explode. These results suggest that the random environment can neither cause the explosion nor prevent the explosion  from happening.

The rest of the paper is arranged as follows. We introduce the CBLE with competition   and the Foster-Lyapunov type criteria for the explosion and nonexplosion of
CBLE with competition in Section 2. Our main results are stated and proved in Section 3.


\section{CBLEs with competition and the Foster-Lyapunov type criteria}
\setcounter{equation}{0}

Let $(\Omega,\F, (\F_t)_{t\geq 0},\P)$ be a filtered probability space satisfying the usual hypotheses.
Let $\phi$ be a branching mechanism given by
\beqlb\label{phi1}
\phi(\lz)=-b_1\lz+ b_2^2\lz^2+\int_0^{\infty}(\e^{-\lz z}-1+\lz z\textbf{1}_{\{z< 1\}})\mu(dz), \quad \lz\geq 0,
\eeqlb
where $b_1, b_2\in\R$ and $(1\wedge z^2)\mu(dz)$ is a finite measure on $(0,\infty)$. To model the mechanism of random environment,
let $(L(t))_{t\geq 0}$ be a L\'{e}vy process with L\'{e}vy-It\^{o} decomposition:
 \beqlb\label{levy1}
L(t)\ar=\ar\beta t+\sigma B^{(e)}(t)+\int_0^t\int_{[-1,1]}(e^z-1)\tN^{(e)}(ds,dz)\cr \ar\ar+\int_0^t\int_{[-1,1]^c}(e^z-1)N^{(e)}(ds,dz),
\eeqlb
where $\beta\in \R, \sigma\geq 0$,  $(B^{(e)}(t))_{t\geq 0}$ is a Brownian motion,
$N^{(e)}(ds,dz)$ is a Poisson random measure on $\R_+\times\R$ with intensity $ds\nu(dz)$ satisfying
$\int_{\R}(1\wedge z^2)\nu(dz)<\iz $ and $\tN^{(e)}(ds,dz)=N^{(e)}(ds,dz)-ds\nu(dz)$.

Let $b_0(y)$ be a competition mechanism, that is, $y\mapsto b_0(y)$ is a continuous non-decreasing function on $[0,\infty)$ with
$b_0(0)=0$.
A CBLE with competition can be constructed as the unique strong solution of  following stochastic equation:
\beqlb\label{EQ1.1a}
Y_t\ar=\ar Y_0+\int_0^t(b_1Y_s-b_0(Y_s))ds+\int_0^t\sqrt{2b_2^2Y_s}dB^{(b)}(s)  +\int_0^t\int_0^{1}\int_0^{Y_{s-}}z\tN^{(b)}(ds,dz,du)\cr
\ar\ar + \int_0^t\int_1^{\iz}\int_0^{Y_{s-}}zN^{(b)}(ds,dz,du)+\int_0^tY_{s-}dL(s),
\eeqlb
where  $(B^{(b)}(t))_{t\geq 0}$ is a Brownian motion,
 $N^{(b)}(ds,dz,du)$ is a Poisson random measure on $\R_+^3$ with intensity $ds\mu(dz)du$,
and $\tN^{(b)}(ds,dz,du)=N^{(b)}(ds,dz,du)-ds\mu(dz)du$.
We also assume that $(B^{(b)}(t))_{t\geq 0}$, $(B^{(e)}(t))_{t\geq 0}$, $N^{(b)}(ds,dz,du)$ and $N^{(e)}(ds,dz)$ are independent of each other.

For $u\geq 0$, let
\beqnn
\tau_u^-:=\inf\{t\geq 0:Y(t)\leq u\}\quad\mbox{and} \quad \tau_u^+:=\inf\{t\geq 0:Y(t)\geq u\}
\eeqnn
and
\beqnn
\tau_{0}:=\tau_0^-\quad\mbox{and} \quad \tau_{\iz}:=\lim_{u\to\infty}\tau^+_u
\eeqnn
with the convention $\inf \emptyset =\iz$.
Throughout this paper, we use notation
\beqnn
\P_{y_0}\{\, \cdot\, \}=\P\{\,  \cdot \,  |Y_0=y_0\} \quad \mbox{and}\quad \E_{y_0}[\,  \cdot \, ]=\E[\,  \cdot \,  |Y_0=y_0], \qquad y_0\geq 0.
\eeqnn

 A $[0, \infty]$-valued process $(Y_t)_{t\geq 0}$ with c\`{a}dl\`{a}g path is a solution to SDE (\ref{EQ1.1a}) if it satisfies  (\ref{EQ1.1a}) up  to explosion time $\tau_\infty$ and $Y_t:=\infty$ for all $t\geq \tau_\infty$. It is known that SDE (\ref{EQ1.1a}) has a unique non-negative strong solution; see Theorem 1 of \cite{PP18}.

Let $L$ be the generator of the process $(Y_t)_{t\geq 0}$.
By It\^{o}'s formula we get for $g\in C^2(\mathbb{R})$,
 \beqlb\label{L3}
Lg(y)\ar=\ar [\beta y+ b_1y-b_0(y)]g'(y)+(\frac{1}{2}\sigma^2 y^2+b_2^2y)g''(y)\cr
\ar\ar + y\int_0^{1}[g(y+z)-g(y)-g'(y)z]\mu(dz) + y\int_1^{\iz}[g(y+z)-g(y)]\mu(dz)\cr
\ar\ar +\int_{[-1,1]}[g(y\e^z)-g(y)-y(\e^z-1)g'(y)]\nu(dz)\cr
\ar\ar +\int_{[-1,1]^c}[g(y\e^z)-g(y)]\nu(dz).
\eeqlb

The Foster-Lyapunov criteria are first used to classify the boundaries for Markov chains via conditions on the generators; see Meyn and Tweedie \cite{MT93} and Chen \cite{C04} for earlier results. These techniques are applied in \cite{LYZ19, MYZ21} and Ren et al. \cite{RXYZ21} to study the boundary behaviours for SDEs associated to continuous-state branching processes.
By a simple modification of the proof of Proposition 2.1  in \cite{RXYZ21} we have the following proposition on solution $Y$ to SDE (\ref{EQ1.1a}):
\bprop\label{p2.3}
If there exist a sequence of strictly positive constants $(d_{n})_{n\geq 1}$
and non-negative functions  $g_n\in C^2((0,\iz))$ satisfying, for all large enough $n\geq 1$,
\begin{enumerate}[(i)]
  \item $\lim_{y\to\infty}g_n(y)=\infty$,
  \item $Lg_n(y)\leq d_{n}g_n(y)$ for all $y\in[1/n,\iz)$,
\end{enumerate}
then $\P_{y_0}\{\tau_{\iz}<\iz\}=0$ for any $y_0>0$.
\eprop

\bprop\label{p2.4}
If  there exist a non-negative bounded and strictly increasing function $g\in C^2((0,\iz))$ and
 positive constants $d_0, \by>0$ satisfying
$$Lg(y)\geq d_0g(y)\quad\text{ for all}\quad y\geq \by,$$
then $\P_{y_0}\{\tau_{\iz}<\iz\}>0$ for any $y_0>\by$.
\eprop

\proof
By the proof of Proposition 2.2 in Ren et al.\cite{RXYZ21} we have, for any $m>\by$,
\beqnn
t\mapsto M_t\ar:=\ar g(Y_{t\wz\tau_m^+\wz\tau_{\by}^-}) \e^{-d_0t} +\int_0^tg(Y_{s\wz\tau_m^+\wz\tau_{\by}^-})d_0\e^{-d_0s}ds\cr \ar\ar-\int_{0}^{t}\e^{-d_0s}Lg(Y_{s})1_{\{s\leq\tau_m^+\wz\tau_{\by}^-\}} ds
\eeqnn
is a martingale. Then for any $m>y_0>\by$,
\beqnn
\ar\ar\E_{y_0}\Big[ g(Y_{t\wz\tau_m^+\wz\tau_{\by}^-})\e^{-d_0t} \Big]
+\int_{0}^{t}\E_{y_0}\Big[d_0\e^{-d_0s}g(Y_{s\wz\tau_m^+\wz\tau_{\by}^-})\Big]ds\cr
\ar\ar=g(y_0)+\int_{0}^{t}\E_{y_0}\Big[\e^{-d_0s}Lg(Y_{s})1_{\{s\leq\tau_m^+\wz\tau_{\by}^-\}} \Big]ds.
\eeqnn
Letting $t\to\iz$, by the assumptions and the dominated convergence theorem  we have
\beqnn
 \int_{0}^{\iz}\E_{y_0}\Big[d_0\e^{-d_0s}g(Y_{s\wz\tau_m^+\wz\tau_{\by}^-})\Big]ds
 \ar=\ar g(y_0)+\int_{0}^{\iz}\E_{y_0}\Big[\e^{-d_0s}Lg(Y_{s})1_{\{s\leq\tau_m^+\wz\tau_{\by}^-\}} \Big]ds\cr
 \ar\geq \ar g(y_0)+\int_{0}^{\iz}\E_{y_0}\Big[\e^{-d_0s}d_0g(Y_{s})1_{\{s\leq\tau_m^+\wz\tau_{\by}^-\}} \Big]ds,
\eeqnn
which implies
\beqlb\label{ineq1}
g(y_0)\ar\leq\ar\E_{y_0}\Big[\int_{\tau_m^+\wz\tau_{\by}^-}^{\iz}d_0\e^{-d_0s}g(Y_{\tau_m^+\wz\tau_{\by}^-})ds\Big]\cr
\ar=\ar\E_{y_0}\Big[g(Y_{\tau_m^+\wz\tau_{\by}^-})\e^{-d_0(\tau_m^+\wz\tau_{\by}^-)}\Big] \cr
\ar=\ar \E_{y_0}\Big[g(Y_{\tau_m^+})\e^{-d_0\tau_m^+}1_{\{\tau_m^+<\tau_{\by}^-\}}\Big]
+\E_{y_0}\Big[g(Y_{\tau_{\by}^-})\e^{-d_0\tau_{\by}^-}1_{\{\tau_m^+>\tau_{\by}^-\}}\Big].
\eeqlb
Since $t\mapsto Y_t$ is right continuous, then $Y_{\tau_{\by}^-}\leq \by<y_0<m\leq Y_{\tau_m^+}$.
Notice that $g$ is non-negative bounded and strictly increasing. Then
\beqnn
g(y_0)\leq \bar{g}\E_{y_0}\Big[1_{\{\tau_m^+<\tau_{\by}^-\}}\e^{-d_0{\tau_m^+}}\Big]+g(\by).
\eeqnn
where $\bar{g}:=\sup_{y}g(y)<\infty$.  Letting $m\to\iz$, we get
\beqnn
g(y_0)\leq \bar{g}\E_{y_0}\Big[1_{\{\tau_{\iz}\leq\tau_{\by}^-\}}\e^{-d_0{\tau_{\iz}}}\Big]+g(\by).
\eeqnn
That is
\beqlb\label{kk}
\bar{g}\E_{y_0}\Big[1_{\{\tau_{\iz}\leq\tau_{\by}^-\}}1_{\{\tau_{\iz}<\iz\}}\Big]\geq g(y_0)-g(\by),
\eeqlb
which implies $$\P_{y_0}\{\tau_{\iz}<\iz\}\geq\frac{g(y_0)-g(\by)}{\bar{g}}>0.$$
This proves the desired result.\qed

\section{Main results}
\setcounter{equation}{0}

In this Section we provide the sufficient conditions for explosion and non-explosion of the CBLE  with competition.
 Let $(Y_t)_{t\geq 0}$  be the unique strong solution of (\ref{EQ1.1a}).

Let $B(p,q)$ denote the Beta function with parameters $p,q>0$.
By integration by parts and L'H\^{o}pital's rule it is not hard to see the following:
 \blem\label{l3.1}
  For any $\dz,y>0$ and $\az\in(0,1)$, we have
        \beqlb\label{beta1}
          \int_0^{\iz}[(y+z)^{-\dz}-y^{-\dz}]z^{-1-\az}dz=-\dz c_{\az,\dz} y^{-\az-\dz}
        \eeqlb
         and
        \beqlb\label{beta1ln}
         \int_0^{\iz}[\ln(y+z)-\ln y]z^{-1-\az}dz=c_{\az,0} y^{-\az},
        \eeqlb
       where $c_{\az,\dz}:={\az}^{-1} B(\az+\dz,1-\az)$.
 \elem

\begin{remark}      Note that $c_{\az,0}=\frac{\pi}{\az\sin(\az\pi)}$.
\end{remark}

For two $\sigma$-finite measures $\mu_1$ and $\mu_2$ on $(0,\infty)$,
we write $\mu_1(dz)\leq \mu_2(dz)$ if $\mu_1(B)\leq \mu_2(B)$ for any Borel set  $B$ in $(0,\infty)$. We first present a sufficient condition on explosion of the solution $Y$ to SDE (\ref{EQ1.1a}).

\bthm\label{t3.1} Suppose that there exist constants $b_0\geq 0$, $q_0\in\R$, $\bar{a}, A>0$  and $\az\in(0,1)$
such that
\[b_0(y)\leq b_0y^{q_0}\quad\text{  for all}\quad y\geq A\quad\text{ and}\quad
\bar{a}z^{-1-\az}\emph{\textbf{1}}_{\{z\geq A\}}dz\leq\emph{\textbf{1}}_{\{z\geq A\}}\mu(dz).\]
Then  $\P_{y_0}\{\tau_{\iz}<\iz\}>0$ for large enough $y_0>0$   if one of the following conditions holds:
\begin{enumerate}[(i)]
  \item $b_0=0$;
  \item  $q_0<2-\az$  and  $b_0>0$;
  \item $q_0=2-\az$ and $0<b_0<\bar{a}c_{\az,0}$.
\end{enumerate}
\ethm

\proof
Without loss of generality we can assume that $A>1$.
Given $\dz\in(0,\infty)$, let $g(y)=\e^{-y^{-\dz}}$ for $y\geq 0$.
Then
\beqnn
g'(y)=\dz y^{-\dz-1}g(y) \quad\mbox{and}\quad g''(y)=[\dz^2 y^{-2\dz-2}-\dz(1+\dz)y^{-\dz-2}]g(y).
\eeqnn
It follows that  $g'(y)>0$ and $g''(y)> -\dz(1+\dz)y^{-\dz-2}g(y)$, which implies that
\beqlb\label{3.1}
\int_0^{1}[g(y+z)-g(y)-zg'(y)]\mu(dz)\ar=\ar\int_0^{1}z^2g''(\xi_1)\mu(dz)\cr
\ar \geq\ar -\dz(1+\dz)y^{-\dz-2}g(y+1)\int_0^{1}z^2\mu(dz)
\eeqlb
for some $\xi_1\in[y,y+1]$, and
\beqlb\label{3.2}
\ar\ar \int_1^{A}[g(y+z)-g(y)]\mu(dz)\geq 0.
\eeqlb
Moreover, by the assumptions  and (\ref{beta1}) we have
\beqlb\label{3.3}
\ar\ar \int_A^{\iz}[g(y+z)-g(y)]\mu(dz)\cr
\ar\ar \geq  \bar{a}\int_A^{\iz}[g(y+z)-g(y)] z^{-1-\az}dz\cr
\ar\ar = \bar{a}g(y)\int_A^{\iz}\Big[\e^{-(y+z)^{-\dz}+y^{-\dz}}-1\Big]z^{-1-\az}dz\cr
\ar\ar \geq \bar{a} g(y)\int_A^{\iz}[-(y+z)^{-\dz}+y^{-\dz}]z^{-1-\az}dz\cr
\ar\ar = \bar{a} g(y)\int_0^{\iz}[-(y+z)^{-\dz}+y^{-\dz}] z^{-1-\az}dz- \bar{a} g(y)\int_0^{A}[-(y+z)^{-\dz}+y^{-\dz}] z^{-1-\az}dz\cr
\ar\ar = \bar{a} g(y)\dz c_{\az,\dz}y^{-\az-\dz}-\bar{a} g(y)\dz\xi_2^{-1-\dz}\int_0^{A}z^{-\az}dz\cr
\ar\ar \geq g(y)\dz \bar{a}c_{\az,\dz}y^{-\az-\dz}-g(y)\dz\bar{a}(1-\az)^{-1}A^{1-\az}y^{-1-\dz},
\eeqlb
where $\xi_2\in[y,y+A]$.
In view of (\ref{3.1})-(\ref{3.3}) we get
\beqlb\label{3.4}
\ar\ar y\int_0^{1}[g(y+z)-g(y)-zg'(y)]\mu(dz)+y\int_1^{\iz}[g(y+z)-g(y)]\mu(dz)\cr
\ar\ar \geq g(y)\dz \bar{a}c_{\az,\dz}y^{1-\az-\dz}-g(y)\dz\bar{a}(1-\az)^{-1}A^{1-\az}y^{-\dz}\cr
\ar\ar \quad -g(y+1)\dz(1+\dz)y^{-\dz-1}\int_0^{1}z^2\mu(dz).
\eeqlb

On the other hand, it is obvious that for any fixed $\dz>0$
there exists a large enough $y_{\dz}>0$ such that $g''(y) <0$ for all $y>y_{\dz}$.
Since $|\e^z-1|\leq 3|z|$ for $z\in [-1,1]$, then for all $y>\e y_{\dz}$
\beqnn
\ar\ar\int_{[-1,1]}[g(y\e^z)-g(y)-y(\e^z-1)g'(y)]\nu(dz)\cr
\ar\ar \geq 9 y^2[g''(\xi_3)\int_{-1}^0z^2\nu(dz)+g''(\xi_4)\int_0^1z^2\nu(dz)]
\eeqnn
for some $\xi_3\in[y\e^{-1}, y]$ and $\xi_4\in[y, y\e]$. This together with $g''(y)>-\dz(1+\dz)y^{-\dz-2}$ yields, for all $y>\e y_{\dz}$,
\beqlb\label{3.5}
\ar\ar\int_{[-1,1]}[g(y\e^z)-g(y)-y(\e^z-1)g'(y)]\nu(dz)\cr
\ar\ar \geq -9\dz(1+\dz)y^2\Big[\xi_3^{-\dz-2}\int_{-1}^0z^2\nu(dz)+\xi_4^{-\dz-2}\int_0^1z^2\nu(dz)\Big]\cr
\ar\ar \geq -9\dz(1+\dz)y^{-\dz}\Big[e^{\dz+2}\int_{-1}^0z^2\nu(dz)+\int_0^1z^2\nu(dz)\Big].
\eeqlb
Moreover, since $g$ is strictly increasing and takes values in $[0,1]$,
we have
$$\int_1^{\infty}[g(y\e^z)-g(y)]\nu(dz)\geq 0.$$
Indeed,
\beqnn
\int_1^{\infty}[g(y\e^z)-g(y)]\nu(dz)\ar=\ar\int_1^{\infty}[\e^{-(y\e^z)^{-\dz}}-\e^{-y^{-\dz}}]\nu(dz)\cr
\ar=\ar g(y)\int_1^{\infty}[\e^{-(y^{-\dz}\e^{-\dz z})+y^{-\dz}}-1]\nu(dz)\cr
\ar=\ar g(y)\int_1^{\infty}[\e^{y^{-\dz}(1-\e^{-\dz z})}-1]\nu(dz)\cr
\ar\to\ar 0,\quad \mbox{as}~ y\to\infty.
\eeqnn
It follows that
\beqlb\label{3.6}
\int_{[-1,1]^c}[g(y\e^z)-g(y)]\nu(dz)\ar\geq\ar \int_{-\infty}^{-1}[g(y\e^z)-g(y)]\nu(dz)\geq -g(y)\nu((-\infty,-1]).
\eeqlb
Combining  (\ref{L3}) and (\ref{3.4})-(\ref{3.6}),  we have, for all $y$ large enough,
\beqlb\label{ineqL}
Lg(y)\ar=\ar [\beta y+ b_1y-b_0(y)]\dz y^{-\dz-1}g(y)
+(\frac{1}{2}\sigma^2 y^2+b_2^2y)[\dz^2 y^{-2\dz-2}-\dz(1+\dz)y^{-\dz-2}]g(y)\cr
\ar\ar +y\int_0^{1}[g(y+z)-g(y)-g'(y)z]\mu(dz) +y\int_1^{\iz}[g(y+z)-g(y)]\mu(dz)\cr
\ar\ar +\int_{[-1,1]}[g(y\e^z)-g(y)-y(\e^z-1)g'(y)]\nu(dz)+\int_{[-1,1]^c}[g(y\e^z)-g(y)]\nu(dz)\cr
\ar\geq\ar g(y)\dz \Big[(\beta+b_1) y^{-\dz}-b_0y^{q_0-1-\dz}-\frac{1}{2}\sigma^2 (1+\dz)y^{-\dz}-b_2^2(1+\dz)y^{-1-\dz}\cr
\ar\ar +\bar{a}c_{\az,\dz}y^{1-\az-\dz} -\bar{a}(1-\az)^{-1}A^{1-\az}y^{-\dz}\cr
\ar\ar -g(y)^{-1}g(y+1)(1+\dz)y^{-\dz-1}\int_0^{1}z^2\mu(dz)\cr
\ar\ar -\dz^{-1}\nu((-\infty,-1])-9(1+\dz)y^{-\dz}e^{y^{-\dz}}e^{\dz+2}\int_{-1}^1z^2\nu(dz)\Big]\cr
\ar=\ar g(y)\dz[\bar{a}c_{\az,\dz}y^{1-\az-\dz}-b_0y^{q_0-1-\dz}-\dz^{-1}\nu((-\infty,-1])-O(y^{-\dz})]\cr
\ar=:\ar g(y)\dz G_{\dz}(y),
\eeqlb
where $O(y^{-\dz})\to 0$ as $y\to\iz$ for any $\dz>0$.

Since $\az<1$, we can first choose $\dz$ small enough such that $1-\az-\dz>0$. If condition (ii)  holds, then $1-\az-\dz>q_0-1-\dz$. Therefore, $G_{\dz}(y)\to \infty$ as $y\to \infty$ under condition (i) or (ii).
If condition (iii) holds, we can choose $\dz$ small enough
such that $1-\az-\dz>0$ and $b_0<\bar{a}c_{\az,\dz}$, then we also have $G_{\dz}(y)\to \infty$ as $y\to \infty$.
This together with (\ref{ineqL}) implies that there is a $\by$ large enough such that $Lg(y)\geq g(y)$ for all $y\geq \by$.
By Proposition \ref{p2.4} we obtain the desired result.
\qed

We next present a sufficient condition   on nonexplosion of process $Y$.

\bthm\label{t3.2} Suppose that there exist constants  $b_0\geq 0$, $q_0\in\R$, $\bar{a},A>0$ and $\az\in(0,2)$
such that
\[b_0(y)\emph{\textbf{1}}_{\{\az< 1\}}\geq b_0y^{q_0}\emph{\textbf{1}}_{\{\az< 1\}}\quad\text{ for all}\quad y\geq A\quad\text{ and}\quad
\emph{\textbf{1}}_{\{z\geq A\}}\mu(dz)\leq\bar{a}z^{-1-\az}\emph{\textbf{1}}_{\{z\geq A\}}dz.\]
Then
$\P_{y_0}\{\tau_{\iz}<\iz\}=0$ for any $y_0>0$   if one of the following conditions holds:
\begin{enumerate}[(i)]
  \item  $\az\geq 1$;
  \item  $q_0>2-\az>1$  and $b_0>0$;
  \item $q_0=2-\az>1$ and $b_0\geq\bar{a}c_{\az,0}$.
\end{enumerate}
\ethm

\proof
For $k\geq 2$, we consider the following stochastic equation:
\beqlb\label{EQ1.2a}
Y_t^{(k)}\ar=\ar Y_0^{(k)}+\int_0^t(\beta Y_s^{(k)}+b_1Y_s^{(k)}-b_0(Y_s^{(k)}))ds
+\int_0^t\sqrt{2b_2^2Y_s^{(k)}}dB^{(b)}(s) \cr
\ar\ar+\int_0^t\sigma Y_s^{(k)}dB^{(e)}(s)
+\int_0^t\int_0^{1}\int_0^{Y_{s-}^{(k)}}z\tN^{(b)}(ds,dz,du)\cr
\ar\ar + \int_0^t\int_1^{\iz}\int_0^{Y_{s-}^{(k)}}zN^{(b)}(ds,dz,du)
+\int_0^t\int_{[-1,1]}Y_{s-}^{(k)}(\e^z-1)\tN^{(e)}(ds,dz) \cr
\ar\ar + \int_0^t\int_{(-\infty,-1)\cup(1,k]}Y_{s-}^{(k)}(\e^z-1)N^{(e)}(ds,dz).
\eeqlb
By Theorem 1 in Palau and Pardo \cite{PP18}, for any $k\geq 2$, equation (\ref{EQ1.2a}) has a unique strong solution $(Y_t^{(k)})_{t\geq 0}$.
 Clearly, $(Y_t^{(k)})_{t\geq 0}$ consists in truncation of large jumps due to environment.
Let $L_k$ be the generator of $(Y_t^{(k)})_{t\geq 0}$. Then
\beqlb\label{L3a}
L_kg(y)\ar=\ar (\beta y+ b_1y-b_0(y))g'(y)+(\frac{1}{2}\sigma^2 y^2+b_2^2y)g''(y)\cr
\ar\ar + y\int_0^{1}[g(y+z)-g(y)-g'(y)z]\mu(dz) + y\int_1^{\iz}[g(y+z)-g(y)]\mu(dz)\cr
\ar\ar +\int_{[-1,1]}[g(y\e^z)-g(y)-y(\e^z-1)g'(y)]\nu(dz)\cr
\ar\ar +\int_{(-\infty,-1)\cup(1,k]}[g(y\e^z)-g(y)]\nu(dz).
\eeqlb

We first prove that for any fixed $k\geq 2$, process $(Y_t^{(k)})_{t\geq 0}$ does not explode.
Without loss of generality we assume $A>1$.

 For $n\geq 9$, let $g_n\in C^2((0,\iz))$ be a non-decreasing function with
$g_n(y)=\ln\ln (n^2y)$ for $y\geq 1/(n\e)$ and $g_n(y)=0$ for $y\leq 1/(2ne)$.
Then for any $y\geq 1/n$,
\beqnn
g_n'(y)=(\ln n^2y)^{-1}y^{-1}>0\quad\text{and}\quad g_n''(y)=-(\ln n^2y)^{-2}y^{-2}-(\ln n^2y)^{-1}y^{-2}<0.
\eeqnn
It follows that
\beqlb\label{3.7}
\ar\ar\int_0^{1}[g_n(y+z)-g_n(y)-zg_n'(y)]\mu(dz)\leq 0
\eeqlb
and
\beqlb\label{3.8}
\int_1^{A}[g_n(y+z)-g_n(y)]\mu(dz)\ar\leq\ar  y^{-1}(\ln n^2y)^{-1}\int_1^{A}z\mu(dz)
\eeqlb
for $y\geq 1/n$.
By the assumption on $\mu$ we have
\beqlb\label{3.9}
\int_A^{\iz}[g_n(y+z)-g_n(y)] \mu(dz)\ar\leq\ar \bar{a}\int_A^{\iz}[g_n(y+z)-g_n(y)] z^{-1-\az}dz.
\eeqlb
If $\az\geq 1$, by integration by parts and L'H\^{o}pital's rule we get for $y\geq 1/n$,
\beqlb\label{3.10}
\ar\ar\int_A^{\iz}[g_n(y+z)-g_n(y)] z^{-1-\az}dz\cr
\ar\ar\leq\int_A^{\iz}[g_n(y+z)-g_n(y)] z^{-2}dz\cr
\ar\ar =A^{-1}[\ln\ln n^2(y+A)-\ln\ln n^2y]+\int_A^{\iz}\frac{1}{z(y+z)\ln n^2(y+z)}dz\cr
\ar\ar \leq(y\ln n^2y)^{-1}+(\ln n^2y)^{-1}\int_A^{\iz}\frac{1}{z(y+z)}dz\cr
\ar\ar = y^{-1}(\ln n^2y)^{-1}[1+\ln(1+y/A)].
\eeqlb
From (\ref{3.7})-(\ref{3.10}) we get
\beqlb\label{3.11}
\ar\ar y\int_0^{1}[g_n(y+z)-g_n(y)-zg_n'(y)]\mu(dz)+ y\int_1^{\iz}[g_n(y+z)-g_n(y)]\mu(dz)\cr
\ar\ar\leq y\int_1^{A}[g_n(y+z)-g_n(y)]\mu(dz)+y\int_A^{\iz}[g_n(y+z)-g_n(y)]\mu(dz)\cr
\ar\ar\leq (\ln n^2y)^{-1} \Big[\int_1^{A}z\mu(dz)+\bar{a}(1+\ln(1+y/A))\Big].
\eeqlb
On the other hand,  since $g_n''(y)\leq 0$ for $y\geq 1/(n\e)$,
\beqlb\label{3.12}
\int_{[-1,1]}[g_n(y\e^z)-g_n(y)-y(\e^z-1)g_n'(y)]\nu(dz)\leq 0
\eeqlb
for $y\geq 1/n$.
Set
\beqnn \gamma_n(y,z):=\ln (\ln  n^2y +z)-\ln (\ln  n^2y)=\ln (1+\frac{z}{\ln  n^2y}).
\eeqnn
Clearly,  $y\mapsto \gamma_n(y,z)$ is strictly decreasing and
$\lim_{y\to\infty} \gamma_n(y,z)=0$ for all $z>0$. Then we can use the monotone convergence to conclude
\beqlb\label{3.13}
\ar\ar\int_{(-\infty,-1)\cup(1,k]}[g_n(y\e^z)-g_n(y)]\nu(dz)\cr
\ar\ar \leq \int_{1}^{k}[\ln (\ln n^2y +z)-\ln (\ln n^2y)]\nu(dz)\cr
\ar\ar =\int_{1}^{k}\gamma_n(y,z)\nu(dz)\to 0, \quad \mbox{as} ~y\to \infty.
\eeqlb

For all $y\geq 1/n$, by (\ref{L3a}) and (\ref{3.11})-(\ref{3.13}) and using $b_0\geq 0$ and $g''<0$ we see that,
if condition (i) holds, then
\beqnn
L_kg_n(y)\ar\leq\ar (\beta y+ b_1y)(\ln n^2y)^{-1}y^{-1}
+(\ln n^2y)^{-1} \Big[\int_1^{A}z\mu(dz)+\bar{a}(1+\ln(1+y/A))\Big]\cr
\ar\ar +\int_{(-\infty,-1)\cup(1,k]}[g_n(y\e^z)-g_n(y)]\nu(dz)\cr
\ar\leq\ar  (\ln n^2y)^{-1}\Big[\beta +b_1+\int_1^{A}z\mu(dz)+\bar{a}(1+\ln(1+y/A))\Big] +\int_{1}^{k}\gamma_n(y,z)\nu(dz)\cr
\ar=:\ar G_{n,k}(y).
\eeqnn
Clearly, for any $n\geq 9$, $G_{n,k}(y)$ converges to some constant as $y\to\infty$ and then $G_{n,k}(y)$ is bounded on $[1/n,\iz)$.
Since $g_n(y)\geq 1$ on $[1/n,\iz)$, then for all $n\geq 9$ there exists a constant $d_n$ such that $L_kg_n(y)\leq d_ng_n(y)$.
By Proposition \ref{p2.3} we have   $(Y_t^{(k)})_{t\geq 0}$ does not explode for all $k\geq 2$.

We now focus on the case that $\az<1$.
Write $a=\ln n^2y$ and $b=\ln n^2(y+z)$. We clearly have $0<a<b$ for $y\geq 1/n$ and then $\ln b-\ln a \leq a^{-1}(b-a)$
 by the concaveness  of the logarithm. Thus,
\beqnn
g_n(y+z)-g_n(y)\leq (\ln n^2y)^{-1}[\ln (y+z)-\ln y], \quad y\geq 1/n.
\eeqnn
This combined with (\ref{beta1ln}) implies
\beqlb\label{3.14}
\int_A^{\iz}[g_n(y+z)-g_n(y)] z^{-1-\az}dz \ar\leq \ar (\ln n^2y)^{-1}\int_0^{\iz}[\ln(y+z)-\ln y]z^{-1-\az}dz \cr
\ar = \ar(\ln n^2y)^{-1}c_{\az,0}y^{-\az}.
\eeqlb
By (\ref{3.7})-(\ref{3.9}) and (\ref{3.14}) we get
\beqlb\label{3.15}
\ar\ar y\int_0^{1}[g_n(y+z)-g_n(y)-zg_n'(y)]\mu(dz)+ y\int_1^{\iz}[g_n(y+z)-g_n(y)]\mu(dz)\cr
\ar\ar\leq (\ln n^2y)^{-1} \Big[\int_1^{A}z\mu(dz)+\bar{a}c_{\az,0}y^{1-\az}\Big].
\eeqlb
For all $y\geq 1/n$, one can use (\ref{L3a}), (\ref{3.12}), (\ref{3.13}) and (\ref{3.15}) to see that
\beqnn
L_kg_n(y)\ar\leq\ar [\beta y+ b_1y-b_0(y)](\ln n^2y)^{-1}y^{-1}
+(\ln n^2y)^{-1} \Big[\int_1^{A}z\mu(dz)+\bar{a}c_{\az,0}y^{1-\az}\Big] \cr
\ar\ar +\int_{(-\infty,-1)\cup(1,k]}[g_n(y\e^z)-g_n(y)]\nu(dz)\cr
\ar\leq\ar  (\ln n^2y)^{-1}\Big[\beta +b_1-b_0(y)y^{-1}+\bar{a}c_{\az,0}y^{1-\az}+ \int_1^{A}z\mu(dz)\Big] +\int_{1}^{k}\gamma_n(y,z)\nu(dz)\cr
\ar=:\ar \bar{G}_{n,k}(y).
\eeqnn
Under the assumption $b_0(y)\geq b_0y^{q_0}$ for all $y\geq A$, if either condition (ii) or condition (iii) holds,
it is not hard to show that for all $k\geq 2$, $y\mapsto\bar{G}_{n,k}(y)$ is bounded above on $[1/n,\infty)$,
and hence $(Y_t^{(k)})_{t\geq 0}$ does not explode by Proposition \ref{p2.3}.

 Now, let $(Y_t)_{t\geq 0}$  be the unique  strong solution of (\ref{EQ1.1a}).  We proceed to show that $(Y_t)_{t\geq 0}$ does not explode.
 Clearly, equation (\ref{EQ1.1a}) can be rewritten as
\beqnn
Y_t\ar=\ar Y_0+\int_0^t(\beta Y_s+b_1Y_s-b_0(Y_s))ds+\int_0^t\sqrt{2b_2^2Y_s}dB^{(b)}(s) +\int_0^t\sigma Y_sdB^{(e)}(s) \cr
\ar\ar+\int_0^t\int_0^{1}\int_0^{Y_{s-}}z\tN^{(b)}(ds,dz,du)+ \int_0^t\int_1^{\iz}\int_0^{Y_{s-}}zN^{(b)}(ds,dz,du)\cr
\ar\ar +\int_0^t\int_{[-1,1]}Y_{s-}(\e^z-1)\tN^{(e)}(ds,dz) + \int_0^t\int_{[-1,1]^c}Y_{s-}(\e^z-1)N^{(e)}(ds,dz).
\eeqnn
Define
\beqnn
 Z(t):=\int_0^t\int_1^{\infty}zN^{(e)}(ds,dz)
\eeqnn
 and
 \beqnn
 \sigma_k:=\inf\{t\geq 0:Z(t)-Z(t-)\geq k\}.
 \eeqnn
Then $\{\sigma_k\}_{k\geq 2}$ is non-decreasing and $\sigma_k\to\infty$ a.s. as $k\to\infty$.
On the other hand, by the definition of $\sigma_k$, it is easy to see that
 $(Y_t)_{t\geq 0}$ satisfies (\ref{EQ1.2a}) on the interval $[0,\sigma_k)$ for all $k\geq 2$.
Then the uniqueness of the solution of (\ref{EQ1.2a}) implies $Y_t=Y_t^{(k)}$ for $t<\sigma_k$.
Since $(Y_t^{(k)})_{t\geq 0}$ does not explode for all $k\geq 2$, $\P_{y_0}\{\tau_{\infty}\geq \sigma_k\}=1$ for all $k\geq 2$,
letting $k\to\infty$ we have $\P_{y_0}\{\tau_{\infty}=\infty\}=1$. That gives the desired result.
\qed

\begin{remark}
It follows from	Theorems \ref{t3.1} and \ref{t3.2} that the L\'evy environment does not seem to be essential for the  explosion to happen or not. Intuitively, this is due to the fact that, in contrast  to the jumps corresponding to the branching mechanism in SDE (\ref{EQ1.1a}),
	 the jumps in the last terms of (\ref{EQ1.1a}) arrive at the same rate as the L\'evy process for the environment and do not speed up when the process $Y$ takes large values.
\end{remark}

In the following corollaries, we consider the special case that
\begin{equation}\label{stable}
\mu(dz)=\bar{a}z^{-1-\az}dz\quad\text{ for constants}\quad \bar{a}>0, \, \az\in (0,2)\quad \text{and for all}\quad z>0.
\end{equation}
Combining Theorems \ref{t3.1} and \ref{t3.2} we  immediately have the following corollaries:
\bcor\label{c3.3} Suppose that (\ref{stable}) holds for $\az\geq 1$. Then    
$\P_{y_0}\{\tau_{\iz}<\iz\}=0$ for all  $y_0>0$.
\ecor

\begin{remark}
Note that the process with $\alpha=1$ corresponds to  Neveu's CBLE with competition whose L\'evy measure $\mu$ for the branching mechanism does not satisfy the finite moment condition,  and the above nonexplosion result is not covered in \cite{BPS21} for the CSLE (without competition).
\end{remark}

\bcor\label{c3.4}  Suppose that (\ref{stable}) holds for $\az< 1$ and there exist constants $q_0\in\R$ and $b_0, A\geq 0$ such that $b_0(y)=b_0y^{q_0}$ for $y\geq A$.
Then $\P_{y_0}\{\tau_{\iz}<\iz\}>0$ for large enough $y_0>0$ if and only if one of the following conditions holds:
\begin{enumerate}[(i)]
  \item  $b_0=0$;
  \item  $q_0<2-\az$ and $b_0>0$;
  \item $q_0=2-\az$ and $0<b_0<\bar{a}c_{\az,0}$.
\end{enumerate}
\ecor

\begin{remark}
	Comparing with the integral test in Theorem 1.2 of \cite{LP21}, in which they only considered the special branching mechanism and Brownian environment  and the logistic competition, that is, $b_0(y)=cy^2$ for some $c\geq 0$, the model we consider is more general  and our results agree with that in  \cite{LP21}.
For example, in the case that $\mu(dz)$ is $\az$-stable with $\az\in(0,1)$ and $b_0(y)=cy^2$,
we can immediately conclude from Corollary \ref{c3.4} that the process does not explode if $c>0$ and the process explodes if $c=0$,	which recovers results for this case by the integral test in  \cite{LP21}.
	
\end{remark}

\noindent\textbf{Acknowledgements}\\
\noindent
Rugang Ma thanks Concordia University where this work was completed during his visit. Xiaowen Zhou thanks Clement Foucart for helpful comments and discussions.


\end{document}